\documentclass[12pt]{article}
\usepackage{full page, latexsym,
amscd,
amsthm,
graphicx, epsfig, amssymb}

\newtheorem{thm}{Theorem}[section]

\newtheorem{prop}[thm]{Proposition}

\newtheorem{rema}[thm]{Remark}

\newcommand{\halmos}{\rule{1ex}{1.4ex}}

\newcommand{\nn}{\nonumber \\}

 \newcommand{\pf}{{\it Proof.}\hspace{2ex}}
 \newcommand{\epf}{\hspace*{\fill}\mbox{$\halmos$}}
 \newcommand{\epfv}{\hspace*{\fill}\mbox{$\halmos$}\vspace{1em}}

\newcommand{\C}{\mathbb{C}}
\newcommand{\Z}{\mathbb{Z}}
\newcommand{\R}{\mathbb{R}}

\newcommand{\N}{\mathbb{N}}

\newcommand{\one}{\mathbf{1}}

\title{ {\bf On the applicability of logarithmic tensor category theory} }
\date{
}
\author{Yi-Zhi Huang}

\begin{document}

\bibliographystyle{alpha}
\maketitle
\begin{abstract}
We give results and observations which allow the application of
the logarithmic tensor category theory of Lepowsky, Zhang and the author 
(\cite{HLZ1}--\cite{HLZ9}) to
more general vertex (operator) algebras and their module categories than those 
studied in a paper by the author (\cite{H2}). 
\end{abstract}

\renewcommand{\theequation}{\thesection.\arabic{equation}}

\section{Introduction}

The logarithmic tensor category theory  of Lepowsky, Zhang and the author 
(\cite{HLZ1}--\cite{HLZ9}) gave a construction of a braided tensor category structure 
(including a ribbon structure) on  a module category $\mathcal{C}$ for a suitable vertex operator
algebra  or more general vertex algebra $V$. In this paper we broaden the applicability of this 
theory.  After referring to the assumptions needed for invoking the theory in \cite{HLZ1}--\cite{HLZ9},
we shall show that these assumptions can be verified in more general settings than before, 
and in addition, we can relax an assumption, thus yielding new families of logarithmic tensor 
categories. 

Before we discuss the precise mathematics, we would like to mention that
the present paper is essentially an addendum to \cite{HLZ1}--\cite{HLZ9} and is readable 
only if the reader consults these papers. The reader will need to consult the specific definitions, results
and, especially,  assumptions
 in \cite{HLZ1}--\cite{HLZ9} discussed in this paper. 

The construction  in \cite{HLZ1}--\cite{HLZ9} uses certain assumptions on 
$V$ and $\mathcal{C}$: Assumptions 10.1 in
\cite{HLZ7} and Assumptions 12.1 
and 12.2 in \cite{HLZ9}. Parts 1 to 5 and most assumptions in Part 7 
of Assumption 10.1 in \cite{HLZ7}, and Assumption 12.1 
in \cite{HLZ9} hold for most of the interesting examples and are relatively easy to verify. 
But there are examples of $V$ and $\mathcal{C}$ for which the first half of Part 6 of Assumption 10.1 in \cite{HLZ7} 
does not hold, that is, the weights of elements of
generalized $V$-modules in $\mathcal{C}$ are not all real and we still would like to 
apply the theory of  \cite{HLZ1}--\cite{HLZ9} to these examples. More importantly, the last part (the
statement that $\mathcal{C}$ be closed under $P(z)$-tensor products for some $z\in \C^{\times}$)
of Part
7 of Assumption 10.1 in \cite{HLZ7}, and Assumption
12.2 in \cite{HLZ9}, are related to some of the 
most important and deep properties of 
$V$ and $\mathcal{C}$. When $V$ and the objects of the category $\mathcal{C}$ 
have trivial $A$-gradings (for $V$) and $\tilde{A}$-grading (for the objects of $\mathcal{C}$),
in other words, when the abelian groups $A$ and $\tilde{A}$ are trivial, 
it was also proved in \cite{HLZ8} that Assumption 12.2 follows from
the condition that every finitely-generated lower bounded doubly-graded generalized
$V$-module
be an object of ${\cal C}$ together with the $C_{1}$-cofiniteness condition 
for objects of $\mathcal{C}$  in the sense of \cite{H1.1}.

In \cite{H2}, the author proved that $V$
and the category $\mathcal{C}$ of grading-restricted 
generalized $V$-modules  (that is, strongly graded generalized $V$-modules with  trivial grading abelian group $\tilde{A}$
in the terminology used in \cite{HLZ1}--\cite{HLZ9})
satisfy all the assumptions in \cite{HLZ7} and \cite{HLZ9}  if
$V$ is a vertex operator algebra (in particular, with trivial $A$) satisfying the following
three conditions: 
\begin{enumerate}

\item $V$ is $C^{a}_{1}$-cofinite in the sense that $V/C^{a}_{1}(V)$ is finite dimensional,
where $C_{1}^{a}(V)$ is the 
subspace of $V$ spanned by the elements of the form 
$u_{n}v$ for $u, v\in V_{+}=\coprod_{n\in \Z_{+}}V_{(n)}$ 
and $L(-1)v$ for $v\in V$. (Here $a$ in the superscript of $C^{a}_{1}$ means ``algebra"
since this is the $C_{1}$-cofiniteness condition for $V$ as a vertex operator algebra, not as a $V$-module.)

\item There exists a positive integer $N$ such that 
$|\Re{(n_{1})}-\Re{(n_{2})}|\le N$ for the lowest
weights $n_{1}$ and $n_{2}$ of any two irreducible $V$-modules
and such that $A_{N}(V)$ is finite dimensional.

\item Every irreducible $V$-module $W$  is $\R$-graded and 
$C_{1}$-cofinite in the sense that $W/C_{1}(W)$ is finite dimensional,
where $C_{1}(W)$ is the 
subspace of $W$ spanned by the elements of the form
$u_{n}w$ for $u\in V_{+}=\coprod_{n\in \Z_{+}}V_{(n)}$ and $w\in W$.

\end{enumerate}
In \cite{H2}, the actions of $L(0)$ on objects of the category are in general not semisimple
and thus intertwining operators among these objects in general are indeed logarithmic. 
If $V$ is of positive energy
(that is, $V_{(n)}=0$ for $n<0$ and $V_{(0)}=\C\one$) and $C_{2}$-cofinite
(that is,  $V/C_{2}(V)$ is finite dimensional
where $C_{2}(V)$ is the 
subspace of $V$ spanned by elements of the form
$u_{-2}v$ for $u\in V_{+}$ and $v\in V$), 
then these three conditions hold.
Since these three conditions, or the  positive energy condition  and $C_{2}$-cofiniteness condition,
are relatively easy to verify and have indeed been verified for many interesting examples of 
vertex operator algebras and their module categories, the results in \cite{H2} 
provide a practical method to apply the logarithmic tensor category theory in \cite{HLZ1}--\cite{HLZ9}.
Up to now, the results in \cite{H2} are still the best general results 
on the verifications of the  main assumptions in \cite{HLZ7} and \cite{HLZ9}.

An initial motivation for our embarking on the creation of the logarithmic tensor category
\cite{HLZ1}--\cite{HLZ9} was to show that the braided tensor categories constructed by 
Kazhdan and Lusztig in \cite{KL1}--\cite{KL5} for suitable module categories for affine 
Lie algebras could be recast and understood as a special case of a new ``logarithmic'' 
generalization of the vertex tensor category theory constructed in \cite{HL0}--\cite{HL3} and 
\cite{H1} for suitable module categories for a vertex operator algebra.  In the new theory, the 
module categories would no longer be completely reducible and the actions of $L(0)$ on the 
modules would no longer be semisimple; in the Kazhdan-Lusztig work, the module categories 
have these properties.  The theory, namely, \cite{HLZ1}--\cite{HLZ9}, 
would be---and in fact is---very general 
and not limited to affine Lie algebras.  In \cite{Z}, Zhang indeed placed the Kazhdan-Lusztig 
construction of braided tensor categories into the setting of vertex operator algebra theory, 
in order to apply an earlier version of \cite{HLZ1}--\cite{HLZ9} to construct these braided tensor categories.  
But Proposition 5.8 in \cite{Z} is wrong.  Fortunately, the mistake is minor and we correct it in the 
present paper. Proposition 5.8 in \cite{Z} is an attempt to verify  
the condition that every finitely-generated lower bounded doubly-graded generalized
$V$-module is an object of the category considered (see the comments above about \cite{HLZ8}
and Assumption 12.2 in \cite{HLZ9}). Also, in  \cite{Z}, the objects in the category are not 
all graded by $\R$. 

It is known that vertex operator algebras associated to, and module categories 
corresponding to, the braided tensor categories constructed by Kazhdan and Lusztig \cite{KL1}--\cite{KL5}
do not satisfy the three 
conditions listed above (mainly the second one, but in some cases also the part of the third condition requiring that 
irreducible modules be $\R$-graded). Also, in recent years,  new examples of interesting vertex (operator) 
algebras and their module categories 
have been constructed and studied. Some of these examples also do not satisfy all the three 
conditions listed above. It is therefore important to generalize the results in \cite{H2},
and even in \cite{HLZ1}--\cite{HLZ9}, so that the logarithmic tensor category theory can be applied
to these cases. 

In this short paper, we give results and observations which can be used to apply 
the logarithmic tensor category theory in
\cite{HLZ1}--\cite{HLZ9} to
more general vertex (operator) algebras and module categories than those 
studied in \cite{H2}. In particular, we correct the mistake in \cite{Z} and make sure that 
the results in \cite{Z} indeed  hold even when the objects of $\mathcal{C}$ are not all graded by $\R$. 
This serves to complete the proof that the braided tensor categories of Kazhdan-Lusztig 
are indeed special cases of the logarithmic tensor category theory of \cite{HLZ1}--\cite{HLZ9}.

\paragraph{Acknowledgments} These remarks are based on the authors' answers to questions
asked by Antun Milas. I am very grateful for him for the questions and discussions. I would also
like to thank James Lepowsky for helpful comments.

\section{A unique expansion set result}

Compared with the study of vertex (operator) algebras, modules and single (logarithmic) intertwining operators, 
the study of products and iterates of two or more (logarithmic) intertwining operators is much more difficult but 
is certainly also much richer and deeper. One needs to prove the convergence of the series obtained from 
these products and iterates and then to use analytic extension and expansion in different regions
to obtained the desired results. The analytic extensions of these products and iterates are in general
multivalued analytic functions for which many of the usual techniques, which work perfectly for the rational functions
used in the study of vertex operator algebras and modules, do not work anymore. Among 
those techniques that do not work for intertwining operators is the Laurent expansion of a single-valued
analytic function defined on an annulus. But we still need to prove that the multivalued analytic functions 
obtained from analytic extensions of the products and iterates of (logarithmic) intertwining operators
can be expanded uniquely as series in powers of the variables and in nonnegative integral powers of 
the logarithms of the variables. In general, such expansions, even if they exist, might not be unique. 
The uniqueness is important for us to construct the (logarithmic) intertwining operators needed in 
our results. For this reason, a notion of unique expansion set was introduced in 
\cite{HLZ6}:

We call a subset $\mathcal{S}$ of $\C\times \C$ a
{\em unique expansion set} if the absolute convergence to $0$ on some
nonempty open subset of ${\mathbb C}^{\times}$ of any series
\[
\sum_{(\alpha,\beta)\in{\cal S}} a_{\alpha,\beta}z^\alpha(\log
z)^\beta, \;\;\; a_{\alpha,\beta} \in {\mathbb C},
\]
where $\log z=\log |z|+i\arg z$ and $0\le \arg z<2\pi$, implies that $a_{\alpha,\beta}=0$ for all $(\alpha,\beta)\in{\cal
S}$. Lemma 14.5 in \cite{H1} can be restated as 
saying that the Cartesian product of a
strictly increasing sequence of real numbers and $\{0\}$ is a unique 
expansion set. Proposition 7.8 in \cite{HLZ6} states that 
for any $N \in \N$, $\mathbb{R}\times\{0,\dots,N\}$ is a unique
expansion set. These results involve only real powers of the variable. 
It is known that $\mathbb{C}\times \{0\}$ is not a unique expansion set.

In this section, we give a generalization of Lemma 14.5 in \cite{H1} to the case that 
the powers of the variable can be complex, with finitely many different imaginary parts. 
This generalization will allow us
to apply the theory in \cite{HLZ1}--\cite{HLZ9}  to suitable
module categories for vertex (operator) algebras whose objects might have 
complex weights. In particular, this generalization justifies those
statements in \cite{Z} that also cover the case that the weights of elements 
of modules might be complex. It is still not known whether Proposition 7.8  in \cite{HLZ6} 
(where the powers of the variables are real) 
can be generalized to the case that the powers of the variable contain finitely many different 
imaginary parts. But for existing examples,
the result below, which generalizes Lemma 14.5 in \cite{H1},  is enough.

\begin{prop}\label{complex-exp-set}
Let $\{n_{i}\}_{i\in \Z_{+}}$ be a
sequence of strictly increasing real numbers $($that is, $n_{i}\in \R$ for $i\in \Z_{+}$ and 
 $n_{1}<n_{2}<n_{3}<\cdots)$ and let
$m_{1},\dots, m_{l}$ be distinct real numbers.  Let 
$$\R^{\{n_{i}\}_{i\in \Z_{+}}}_{m_{1}, \dots, m_{l}}=\{n_{i}+m_{j}\sqrt{-1}
\;|\;i\in \Z_{+}, \; j=1, \dots, l\}.$$ Then for any $N \in \N$, $\R^{\{n_{i}\}_{i\in \Z_{+}}}_{m_{1}, \dots, m_{l}}
\times\{0,\dots,N\}$ is a unique
expansion set. 
\end{prop}
\pf Let $a_{i, j, k}\in \C$ for $i\in \Z_{+}$, $j=1, \dots, l$ and $k=0, \dots, N$, and
suppose that
$$
\sum_{i\in \Z_{+}}\sum_{j=1}^{l}\sum_{k=0}^{N}a_{i, j, k}z^{n_{i}+m_{j}\sqrt{-1}}(\log z)^{k}=
\sum_{i\in \Z_{+}}\sum_{j=1}^{l}\sum_{k=0}^{N}a_{i, j, k}e^{(n_{i}+m_{j}\sqrt{-1})\log z}(\log z)^{k}
$$
is absolutely convergent to $0$ for $z=z_{1}\in \C^{\times}$.  (See Remark \ref{rmk-2.2} below.) 
We want to prove that $a_{i, j,k}=0$ for $i\in\Z_{+}$, $j=1, \dots, l$ and $k=0, \dots, N$. 

We may assume that $m_{1}>\cdots >m_{l}$.  It is sufficient to prove that
$a_{1, j, k}=0$ for $j=1, \dots, l$ and $k=0, \dots, N$, which we proceed to do. 

For $z\in \C$ satisfying $|z|\le|z_{1}|$,
$$\sum_{i\in \Z_{+}}|a_{i, j, k}z^{n_{i}-n_{1}}|=\sum_{i\in \Z_{+}}|a_{i, j, k}||z|^{n_{i}-n_{1}}
\le \sum_{i\in \Z_{+}}|a_{i, j, k}||z|_{1}^{n_{i}-n_{1}}= \sum_{i\in \Z_{+}}|a_{i, j, k}z_{1}^{n_{i}-n_{1}}| $$
 is absolutely and uniformly convergent. Thus for $z\in \C$ satisfying $|z|<|z_{1}|$,
\begin{eqnarray}\label{complex-exp-set-0}
\lefteqn{\sum_{i\in \Z_{+}}\sum_{j=1}^{l}\sum_{k=0}^{N}a_{i, j, k}z^{(n_{i}-n_{1})+(m_{j}-m_{l})\sqrt{-1}}(\log z)^{k-N}}\nn
&&=\sum_{i\in \Z_{+}}\sum_{j=1}^{l}\sum_{k=0}^{N}a_{i, j, k}e^{((n_{i}-n_{1})+(m_{j}-m_{l})\sqrt{-1})\log z}(\log z)^{k-N}
\end{eqnarray}
is absolutely and uniformly convergent to $0$.

Take $\log z=x+y\sqrt{-1}$ where $x\in \R$ and $y\in \R$ such that $ e^{x}< |z_{1}|$
or equivalently $x< \log |z_{1}|$. Then $|z|<|z_{1}|$.
Thus from (\ref{complex-exp-set-0}), for such $x$ and $y$, we have 
\begin{equation}\label{complex-exp-set-1}
\sum_{i\in \Z_{+}}\sum_{j=1}^{l}\sum_{k=0}^{N}a_{i, j, k}e^{(n_{i}-n_{1})x}e^{(n_{i}-n_{1})y\sqrt{-1}}
e^{(m_{j}-m_{l})\sqrt{-1}x}
e^{- (m_{j}-m_{l})y}(x+y\sqrt{-1})^{k-N}=0.
\end{equation}
Let $x$ and $y$ go to $-\infty$ and $\infty$, respectively,  on both sides of (\ref{complex-exp-set-1}).
Since the series in the left-hand side of (\ref{complex-exp-set-1}) is uniformly convergent,
we can take the limit term by term. 
Thus we obtain $a_{1, l, N}
=0.$

Now assume that $a_{1, j, N}=0$ for $j=p+1, \dots, l$. Then 
\begin{eqnarray}\label{complex-exp-set-2}
\lefteqn{\sum_{j=1}^{p}\sum_{k=0}^{N}a_{1, j, k}
e^{(m_{j}-m_{p})\sqrt{-1}x}
e^{- (m_{j}-m_{p})y}(x+y\sqrt{-1})^{k-N}}\nn
&&\quad +\sum_{i\in \Z_{+}+1}\sum_{j=1}^{l}\sum_{k=0}^{N}a_{i, j, k}e^{(n_{i}-n_{1})x}e^{(n_{i}-n_{1})y\sqrt{-1}}
e^{(m_{j}-m_{p})\sqrt{-1}x}
e^{- (m_{j}-m_{p})y}(x+y\sqrt{-1})^{k-N}\nn
&&=0.
\end{eqnarray}
Let $y=\log |x|$ in (\ref{complex-exp-set-2}) and then let $x$ go to $-\infty$  on both sides of (\ref{complex-exp-set-2}).
Again by taking the limit on the left-hand side term by term,
we obtain $a_{1, p, N}=0$. Thus we have $a_{1, j, N}=0$ for $j=1, \dots, l$.

It follows that $a_{1, j, k}=0$ for $j=1, \dots, l$ and $k=0, \dots, N-1$ as well.
\epfv

\begin{rema}\label{rmk-2.2}
{\rm From the proof of Proposition \ref{complex-exp-set}, we see that we proved a stronger result.
In fact, in the proof, we have assumed that  the series is absolutely convergent only 
for a particular number  $z=z_{1}\in \C^{\times}$, not for $z$ in a nonempty open subset of $\C^{\times}$. This assumption 
is weaker than the set to be a unique expansion set. Thus we have proved a result stronger than the 
statement of Proposition \ref{complex-exp-set}.}
\end{rema}

\section{A sufficient condition for the expansion condition}

In the proof of Theorem 4.12 in \cite{H2}, the author used the three conditions 
(in fact only the first two conditions) listed in the introduction 
to prove that  every finitely-generated lower-truncated generalized $V$-module
is in the category of grading-restricted generalized $V$-modules (where the grading 
abelian groups $A$ and $\tilde{A}$ are trivial). The convergence and extension property 
for products and iterates of logarithmic intertwining operators also holds by Theorem 11.8 in \cite{HLZ8}. 
These verify the two conditions needed in Theorem 11.4 in \cite{HLZ8}. 
Then by Theorem 11.4 and Theorem 11.8  in \cite{HLZ8}, Assumption 12.2 in \cite{HLZ9} holds.

When the first two conditions listed in the introduction do not hold, especially when the 
second condition does not hold, Condition 1 in Theorem 11.4 in \cite{HLZ8}  might not hold
and thus we cannot use this theorem. 
But when Assumption 10.1 in \cite{HLZ7} holds,  Theorem 11.4 in \cite{HLZ8} 
can be easily generalized by observing that 
the proof of Theorem 11.4 in \cite{HLZ8} in fact proves the following stronger
result (stronger in the sense that one assumption is weaker):

\begin{thm}\label{weaker-form-thm-11.4}
Suppose that Assumption 10.1 in \cite{HLZ7} holds and the following two conditions are satisfied:
\begin{enumerate}
\item  For any objects $W_{1}$ and $W_{2}$ of $\mathcal{C}$ and any $z\in \C^{\times}$, if 
the doubly-graded generalized $V$-module $W_{\lambda}$ 
(or a doubly-graded $V$-module when $\mathcal{C}$ is in $\mathcal{M}_{sg}$)
generated by a generalized eigenvector  $\lambda\in 
(W_{1}\otimes W_{2})^{*}$ for 
$L_{P(z)}(0)$ satisfying the $P(z)$-compatibility condition is lower bounded, then $W_{\lambda}$
is an object of $\mathcal{C}$.

\item The convergence and extension property for either products or
iterates holds in ${\cal C}$ (or the convergence and extension
property without logarithms for either products or iterates holds in
${\cal C}$, when ${\cal C}$ is in $\mathcal{M}_{sg}$).
\end{enumerate}
Then the convergence and expansion conditions for intertwining maps in
${\cal C}$ both hold.
\end{thm}
\pf Note that the only place in the proof  of Theorem 11.4 in \cite{HLZ8} 
where the condition that every finitely-generated lower 
bounded generalized $V$-module is in $\mathcal{C}$ (Condition 1 in that theorem)  is used 
 is in the last paragraph 
showing that $W_{\lambda_{n}^{(2)}(w'_{(4)}, w_{(3)})}$ is in $\mathcal{C}$.
But to show that $W_{\lambda_{n}^{(2)}(w'_{(4)}, w_{(3)})}$ is in $\mathcal{C}$,
Condition 1 in the statement of the present theorem is enough. Thus the theorem is proved.
\epf

\begin{rema}\label{another-weaker-form-thm-11.4}
{\rm In Theorem  \ref{weaker-form-thm-11.4}, we assume in particular that
the (generalized) weights are real numbers for any object in $\mathcal{C}$ 
(the first half of Part 6 of Assumption 
10.1 in \cite{HLZ7}). If we replace this part of the assumption  by the assumption
that every object in $\mathcal{C}$ is of finite length,
that is, has a finite composition series in $\mathcal{C}$, then 
because of  Proposition \ref{complex-exp-set}, the conclusion of 
Theorem \ref{weaker-form-thm-11.4} still holds.}
\end{rema}

Theorem \ref{weaker-form-thm-11.4} and Remark \ref{another-weaker-form-thm-11.4} 
can be used to construct the associativity 
isomorphism when Condition 1 in Theorem 11.4 in \cite{HLZ8}
is not satisfied (in particular,  when $V$ is a vertex operator algebra but
Conditions 1 and 2 in the introduction are not satisfied or, as a special case, 
when $V$ is a vertex operator algebra
that  is not $C_{2}$-cofinite), but  Assumption 10.1 in \cite{HLZ7} 
or the assumptions discussed in Remark \ref{another-weaker-form-thm-11.4}  hold
 and
 the two conditions in Theorem 3.1 are satisfied.

\section{A correction of 
a mistake in a paper of Lin Zhang}

The tensor categories constructed by Kazhdan and Lusztig in \cite{KL} 
correspond to examples of vertex operator 
algebras and module categories that do not satisfy the second condition listed in the introduction
(in particular,  do not satisfy the $C_{2}$-cofiniteness condition).
A construction of these tensor categories of Kazhdan-Lusztig using the logarithmic tensor category 
theory developed in \cite{HLZ1}--\cite{HLZ9} was given in \cite{Z}.

But it was noticed by Milas that one of the propositions in \cite{Z} that is needed in applying 
Theorem 11.4 in \cite{HLZ8} is wrong. Since this proposition in 
\cite{Z} is wrong, we cannot use Theorem 11.4 in \cite{HLZ8}. 
Instead, we use Theorem \ref{weaker-form-thm-11.4} and a result in \cite{KL}
to give an almost trivial correction of the minor 
mistake in \cite{Z}. 

To be more precise, the mistake is that Proposition 5.8  in 
\cite{Z} is wrong. Proposition 5.8 in \cite{Z}, if correct, would verify Condition 1
in Theorem 11.4 in \cite{HLZ8}. To correct this mistake, we need only 
verify Condition 1 in Theorem \ref{weaker-form-thm-11.4} above. We do this by using 
Theorem 7.9  in \cite{KL}. 

\begin{prop}\label{4.1}
For any two objects $W_{1}$ and $W_{2}$ of  $\mathcal{O}_{\kappa}$,
 if the generalized $V_{\hat{\mathfrak{g}}}(\ell, 0)$--module $W_{\lambda}$ 
generated by a generalized eigenvector  $\lambda\in 
(W_{1}\otimes W_{2})^{*}$ for 
$L_{P(z)}(0)$ satisfying the $P(z)$-compatibility condition is lower bounded, then $W_{\lambda}$ is 
an object of $\mathcal{O}_{\kappa}$. 
\end{prop}
\pf
Since $W_{\lambda}$ is lower bounded, by definition, the elements of 
$W_{\lambda}$ must be in $W_{1}\circ_{P(z)}W_{2}$ (see \cite{KL} and \cite{Z}). 
By Theorem 7.9  in \cite{KL}, $W_{1}\circ_{P(z)}W_{2}$ is in $\mathcal{O}_{\kappa}$.
From \cite{KL} (also Theorem 5.1 in \cite{Z}), 
$\mathcal{O}_{\kappa}$ consists of  the $\hat{\mathfrak{g}}$-modules 
 of level $\ell$ having a finite
composition series all of whose irreducible subquotients are of the form
$L(\ell, \mu)$ for various highest weights  $\mu$ for the finite-dimensional 
Lie algebra $\mathfrak{g}$.
Thus $W_{\lambda}$ 
as a submodule of $W_{1}\circ_{P(z)}W_{2}$ is also in $\mathcal{O}_{\kappa}$.
\epf

\begin{rema}\label{4.2}
{\rm  In general, the objects in $\mathcal{O}_{\kappa}$ might have homogeneous subspaces of 
complex weights. In \cite{Z}, it was not justified that the results in \cite{HLZ1}--\cite{HLZ9}
can indeed be applied in this case. What was missing is exactly a unique expansion set result 
that can be used in this case. 
Since any object in $\mathcal{O}_{\kappa}$ is of finite length, for a given object in 
$\mathcal{O}_{\kappa}$, there can only be finitely many different imaginary parts of 
the complex weights of the homogeneous subspaces of the object. Proposition \ref{complex-exp-set}
is exactly the unique expansion set result that we need in this case. Thus Proposition 
\ref{complex-exp-set} fills this minor gap in \cite{Z}. 
}
\end{rema}

Together with \cite{Z}, Proposition \ref{4.1} and Remark \ref{4.2} 
serve to complete the proof that the braided tensor categories of Kazhdan-Lusztig 
are indeed special cases of the logarithmic tensor category theory of \cite{HLZ1}--\cite{HLZ9}.

\noindent {\small \sc Department of Mathematics, Rutgers University,
110 Frelinghuysen Rd., Piscataway, NJ 08854-8019}

\vspace{1em}

\noindent {\em E-mail address}: yzhuang@math.rutgers.edu


\begin{thebibliography}{KWAK2}

\bibitem[H1]{H1} Y.-Z.~Huang, A theory of tensor products for module categories for a 
vertex operator algebra, IV, {\it J. Pure Appl. Alg.} {\bf 100} (1995), 173--216. 

\bibitem[H2]{H1.1} Y.-Z. Huang, Differential equations and
intertwining operators, {\it Comm. Contemp. Math.} {\bf 7} (2005),
375--400.

\bibitem[H3]{H2} Y.-Z.~Huang, Cofiniteness conditions, projective covers and 
the logarithmic tensor product theory, {\it J. Pure Appl. Alg.} {\bf 213} (2009), 458--475.

\bibitem[HL1]{HL0}
Y.-Z. Huang and J. Lepowsky, Toward a
theory of tensor products for representations of a vertex operator
algebra, in: {\em Proc. 20th International Conference on Differential
Geometric Methods in Theoretical Physics, New York, 1991},
ed. S. Catto and A. Rocha, World Scientific, Singapore, 1992, Vol. 1,
344--354.

\bibitem[HL2]{tensorK}
Y.-Z. Huang and J. Lepowsky, Tensor products of modules for a vertex
operator algebras and vertex tensor categories, in:
     {\em Lie Theory and Geometry,
in honor of Bertram Kostant,}
ed. R. Brylinski, J.-L. Brylinski, V. Guillemin, V. Kac,
Birkh\"{a}user, Boston, 1994, 349--383.

\bibitem[HL3]{HL1}
Y.-Z. Huang and J. Lepowsky, A theory of tensor products for module
categories for a vertex operator algebra, I, {\em Selecta Mathematica
(New Series)} {\bf 1} (1995), 699--756.

\bibitem[HL4]{HL2}
Y.-Z. Huang and J. Lepowsky, A theory of tensor products for module
categories for a vertex operator algebra, II, {\em Selecta Mathematica
(New Series)} {\bf 1} (1995), 757--786.

\bibitem[HL5]{HL3}
Y.-Z. Huang and J. Lepowsky, A theory of tensor
products for module categories for a vertex operator algebra, III,
{\em J. Pure Appl. Alg.} {\bf 100} (1995) 141--171.

\bibitem[HLZ1]{HLZ1} Y.-Z.~Huang, J.~Lepowsky and L.~Zhang, A
logarithmic generalization of tensor product theory for modules for a
vertex operator algebra, {\em Internat. J. Math.} {\bf 17} (2006),
975--1012.

\bibitem[HLZ2]{HLZ2}
Y.-Z.~Huang, J.~Lepowsky and L.~Zhang,
Logarithmic tensor category theory for generalized modules for a conformal vertex 
algebra, I: Introduction and strongly graded algebras and their generalized modules,
in: Conformal Field Theories and Tensor Categories, Proceedings of a Workshop Held 
at Beijing International Center for Mathematics Research, ed. C. Bai, J. Fuchs, Y.-Z. 
Huang, L. Kong, I. Runkel and C. Schweigert, Mathematical Lectures from Beijing 
University, Vol. 2, Springer, New York, 2014, 169--248.

\bibitem[HLZ3]{HLZ3} Y.-Z.~Huang, J.~Lepowsky and L.~Zhang, Logarithmic
tensor category theory, II: Logarithmic formal calculus
and properties of logarithmic intertwining operators, to appear; arXiv:1012.4196.

\bibitem[HLZ4]{HLZ4} Y.-Z.~Huang, J.~Lepowsky and L.~Zhang, Logarithmic
tensor category theory, III: Intertwining maps and tensor
product bifunctors, to appear; arXiv:1012.4197.

\bibitem[HLZ5]{HLZ5} Y.-Z.~Huang, J.~Lepowsky and L.~Zhang, Logarithmic
tensor category theory, IV: Constructions of tensor
product bifunctors and the compatibility conditions, to appear; arXiv:1012.4198.

\bibitem[HLZ6]{HLZ6} Y.-Z.~Huang, J.~Lepowsky and L.~Zhang, Logarithmic
tensor category theory, V: Convergence condition for
intertwining maps and the corresponding compatibility
condition, to appear; arXiv:1012.4199.

\bibitem[HLZ7]{HLZ7} Y.-Z.~Huang, J.~Lepowsky and L.~Zhang, Logarithmic
tensor category theory, VI: Expansion condition, associativity of logarithmic
intertwining operators, and the associativity isomorphisms, to appear; arXiv:1012.4202.

\bibitem[HLZ8]{HLZ8} Y.-Z.~Huang, J.~Lepowsky and L.~Zhang, Logarithmic
tensor category theory, VII: Convergence and extension
properties and applications to expansion for intertwining
maps, to appear; arXiv:1110.1929.

\bibitem[HLZ9]{HLZ9} Y.-Z.~Huang, J.~Lepowsky and L.~Zhang, Logarithmic
tensor category theory, VIII: Braided tensor category
structure on categories of generalized modules for a
conformal vertex algebra, to appear; arXiv:1110.1931.

\bibitem[KL1]{KL1}
D. Kazhdan and G. Lusztig,
Affine Lie algebras and quantum groups,
{\it Duke Math. J., IMRN} {\bf 2} (1991), 21--29.

\bibitem[KL2]{KL}
D. Kazhdan and G. Lusztig,
Tensor structures arising {from} affine Lie algebras, I,
{\it J. Amer. Math. Soc.} {\bf 6} (1993), 905--947.

\bibitem[KL3]{KL3}
D. Kazhdan and G. Lusztig,
Tensor structures arising {from} affine Lie algebras, II,
{\it J. Amer. Math. Soc.} {\bf 6} (1993), 949--1011.

\bibitem[KL4]{KL4}
D. Kazhdan and G. Lusztig,
Tensor structures arising {from} affine Lie algebras, III, {\it J.
Amer. Math. Soc.} {\bf 7} (1994), 335--381.

\bibitem[KL5]{KL5}
D. Kazhdan and G. Lusztig,
Tensor structures arising {from} affine Lie algebras, IV,
{\it J. Amer. Math. Soc.} {\bf 7} (1994), 383--453.


\bibitem[Z]{Z} L.~Zhang, Vertex tensor category structure on a
category of Kazhdan-Lusztig, {\it New York J. Math.} {\bf 14} (2008) 261--284.


\end{thebibliography}
\end{document}